\documentclass[10pt]{amsart}
\usepackage{amsmath,amsfonts,amsthm,amstext,amscd}
\usepackage[latin1]{inputenc}
\usepackage[dvips]{graphicx}
\usepackage{array}
\usepackage{hyperref}
\usepackage{newlfont}
\def\rm#1{\mathrm{#1}}
\def\cal#1{\mathcal{#1}}
\def\bb#1{\mathbb{#1}}

\def\lr#1{\left\langle #1 \right\rangle}
\usepackage[all]{xy}

 \newtheorem{theorem}{Theorem}
 \newtheorem{cor}[theorem]{Corollary}
 \newtheorem{lem}[theorem]{Lemma}
 \newtheorem{prop}[theorem]{Proposition}
 \theoremstyle{definition}
 
 \theoremstyle{remark}

 \numberwithin{equation}{section}

\begin{document}

\title[Pulling back the Gromoll-Meyer construction]
 {Pulling back the Gromoll-Meyer construction and models of exotic spheres}

\author{L. D. Speran\c ca}

\address{Departamento de Matem\'atica, UFPR \\
 Setor de Ci\^encias Exatas, Centro Polit\'ecnico, \\
 Caixa Postal 019081,  CEP 81531-990, \\
 Curitiba, PR, Brazil}

\email{lsperanca@ufpr.br}

\thanks{The author was financially supported by FAPESP, grant numbers  2012/25409-6.}


\date{}

\begin{abstract}  
Here we generalize the Gromoll-Meyer construction of an exotic 7-sphere by producing geometric models of exotic 8, 10 and Kervaire spheres as quotients of sphere bundles over spheres by free isometric actions. We give a geometric application at the end.
  \end{abstract}

\maketitle
\section{Introduction}

An {\em exotic sphere} is a differentiable manifold homeomorphic, but not diffeomorphic, to an
standard Euclidean sphere $S^n$; the existence of such manifolds was discovered by Milnor 
\cite{mi}, and in fact he gave a description of (some) 7-dimensional exotic spheres by 
modelling them as 3-sphere bundles over $S^4$ with structure group $SO(4)$. 

A milestone in the study of the geometric structure of exotic spheres was the presentation, 
by Gromoll and Meyer \cite{gm}, of an exotic sphere as a quotient of the Lie 
group $Sp(2)$ of quaternionic 
unitary matrices by an explicit $S^3$ action. More precisely, considering $S^7$ the unitary sphere of the quaternionic plane, $\bb H^2$,
\[Sp(2)=\left\{Q=\begin{pmatrix}a&c\\b&d\end{pmatrix}\in S^7\times S^7~|~\bar ca+\bar db=0\right\}\]

And the group of unit quaternions $S^3$ acts in $Sp(2)$ by
\begin{gather}\label{eq star action}
q\star\begin{pmatrix}a & c\\b & d\end{pmatrix}=\begin{pmatrix}qa\bar q & qc\\qb\bar q & qd\end{pmatrix}
\end{gather}

In fact, the Gromoll-Meyer construction fits inside a richer framework. First observe that 
we have what we call the ``standard'' action of $S^3$ in $Sp(2)$, 
\begin{gather}\label{eq std action}
q\bullet\begin{pmatrix}a & c\\b & d\end{pmatrix}=\begin{pmatrix}a & c\bar q\\b & d\bar q\end{pmatrix}
\end{gather}

These two actions commute, and they give rise to a {\em cross diagram}

\begin{equation}\label{diagram DCD}
\begin{xy}\xymatrix{& S^3\ar@{.}[d]^{\bullet} & \\ S^3\ar@{..}[r]^{\star} &Sp(2)\ar[d]^{pr}\ar[r]^{pr'} &\Sigma^7\\ &S^7&}\end{xy}
\end{equation}

The cross diagram contains a wealth of information about the geometry of the Gromoll-Meyer sphere; 
especially of its geodesics and submanifolds, e.g.\cite{dur,DMR,ADPR,DR}. In particular, note that each one of the actions descends to symmetries of the quotient by the other action.

The main point of this paper is that the Gromoll-Meyer construction can be {\em pullbacked} to give 
explicit quotient models of exotic spheres and related phenomena in dimensions 8 and 10 (Theorem \ref{thm 8 10}). For 7-dimensional 
exotic spheres, it was shown in \cite{DPR} that 
the Gromoll-Meyer construction in facts describes 
{\em all} 7-dimensional exotic spheres, via pullbacking by the Cayley power self-map of $S^7$.

We actually go beyond pullbacking the Gromoll-Meyer construction; in section \ref{sec Kervaire} we give all Kervaire manifolds as quotients by $O(n)$-actions on pullbacks of frame bundle of round spheres.

As an example of geometric application of these presentations we show that the wiedersehen metrics and their relationship with the exotic diffeomorphisms that began in \cite{dur} can also be pullbacked in all of these examples. In \cite{dur}, the geometry of an explicit wiedersehen metric on the Gromoll-Meyer sphere is used to construct a rather symmetric formula for an exotic (degree 1 but not isotopic to the identity) diffeomorphism $\sigma:S^6\to S^6$. If one thinks on homotopy $n$-spheres as manifolds obtained by clutching two discs by a diffeomorphism (as Smale's theorem allows for $n>4$), the group structure given by connected sum of homotopy $n$-spheres is translated as the group structure given by composition in the group of diffeomorphisms of the euclidean $(n-1)$-sphere up to isotopy. One missing desired property of the construction in \cite{dur} is the realization of this group structure as the powers of $\sigma$, in the sense that, the group of homotopy $7$-spheres is the cyclic group with 28 elements although we cannot find an explicity isotopy from $\sigma^{28}$ to the identity map. Here we go further and present a clutching diffeomorphism of the exotic 8-sphere together with an isotopy of its double to the identity (thus realizing the group of homotopy 8-spheres as the cyclic group with two elements).



\section{Pulling Back the Gromoll-Meyer Construction and 8 and 10 dimensional Exotic Spheres}

Consider as in the introduction, the Gromoll-Meyer construction as a pair of commuting free action on $Sp(2)$. We present here a natural way to pull-back both actions simultaneously; we first describe the kind of pullbacks we deal with in general and then specialize to the specific pullbacks that gives 8 and 10 dimensional exotic spheres.

Let's start taking $Sp(2)$ as the principal bundle over $S^7$ defined by action \eqref{eq std action}. Observe that the quotient projection $S^3\stackrel{\bullet}{\cdots} Sp(2)\stackrel{pr}{\to} S^7$ is given by projection into the first column of $Q$, i.e.
\[pr\begin{pmatrix}a & c\\b & d\end{pmatrix}=\begin{pmatrix}a \\b \end{pmatrix}\]

We also recall that, if $f:M\to S^7$ is a map, we can define the \emph{pull-back bundle} $\pi:f^*Sp(2)\to M$ as the bundle with total space
\begin{gather}
f^*Sp(2)=\{(x,Q)\in M\times Sp(2)~|~f(x)=pr(Q)\}
\end{gather}
and projection $\pi(x,Q)=x$. In particular, $(x,Q)\in f^*Sp(2)$ if and only if it is of the form
\begin{equation}\label{eq total pullback}(x,Q)=\left(x,\begin{pmatrix}a(x) & c\\b(x) & d\end{pmatrix}\right) 
\end{equation}
where $a$ and $b$ are comlpetly defined by the identity $(a(x),b(x))^t=\pi(x)\in S^7$. This is a principal bundle with the action 
\[q\bullet \left(x,\begin{pmatrix}a(x) & c\\b(x) & d\end{pmatrix}\right)=\left(x,q\bullet \begin{pmatrix}a(x) & c\\b(x) & d\end{pmatrix}\right)=\left(x,\begin{pmatrix}a(x) & c\bar q\\b(x) & d\bar q\end{pmatrix}\right)\]
which is well-defined since $(c\bar q,d\bar q)$ is quaternionic orthogonal to $(a(x),b(x))$ as it is $(c,d)$. Or, in other words, since $pr(q\bullet Q)=\pi(Q)=f(x)$. Furthermore, if we require some equivariance from $f$, we can indeed pull-back both actions on $Sp(2)$ simultaneously. First observe that, since \eqref{eq star action} commutes with \eqref{eq std action}, it induces the following $S^3$-action on $S^7$, given by projecting $q\star Q$ to the first column:
\begin{gather}\label{eq GM action}
q\cdot \begin{pmatrix}a\\ b\end{pmatrix}=\begin{pmatrix}qa\bar q\\ qb\bar q\end{pmatrix}
\end{gather}
Now, suppose that $M$ is an $S^3$-manifold and that $f:M\to S^7$ is a $S^3$-equivariant map with respect to this action. Then we observe that
\begin{equation}
q\star(x,Q)=\left(qx,\begin{pmatrix}a(qx) & qc\\b(qx) & qd\end{pmatrix}\right)=\left(qx,\begin{pmatrix}qa(x)\bar q & qc\\qb(x)\bar q & qd\end{pmatrix}\right)
\end{equation}
is a well-defined free action on $f^*Sp(2)$ since the rightmost side is clearly a free action and that 
\[f(qx)=q\cdot f(x)=q\cdot pr(Q)=pr(q\star Q).\]

In particular, $f^*Sp(2)$ also fits in a cross-diagram:
\begin{equation}\label{diagram fCD}
\begin{xy}\xymatrix{& S^3\ar@{.}[d]^{\bullet} & \\ S^3\ar@{..}[r]^{\star} &f^*Sp(2)\ar[d]^{\pi}\ar[r]^{\pi'} &M'\\ &M&}\end{xy}
\end{equation}

In what follows we will see how this construction can be used to relate geometry and different constructions of some exotic spheres. 

\subsection{8 and 10 dimensional spheres}

 Consider $S^8$ as the unit sphere on $\bb R\times \bb H^2$ and $f_8:S^8\to S^7$ as
\begin{align}
f_8\begin{pmatrix}\lambda\\x\\w\end{pmatrix}=\frac{1}{\sqrt{\lambda^2+|x|^4+|w|^2}}\begin{pmatrix}\lambda+xi\bar x\\w\end{pmatrix}
\end{align}

This map is equivariant with respect to the action
\begin{equation}
q\cdot\begin{pmatrix}\lambda\\x\\y\end{pmatrix}=\begin{pmatrix}\lambda\\qx\\qw\bar q\end{pmatrix}
\end{equation}
and \eqref{eq GM action}. In particular, as we observed, the manifold $E^{11}=f_8^*Sp(2)$ admits two $S^3$ free action, the usual induced from the principal $\pi^{11}:E^{11}\to S^8$ and a new action 
\begin{eqnarray}
q\star\left(\begin{array}{c}\lambda\\x\\y\end{array}\begin{array}{c} c\\ d\end{array}\right)=\left(\begin{array}{c}\lambda\\qx\\qw\bar q\end{array}\begin{array}{c} qc\\ qd\end{array}\right)
\end{eqnarray}
where, for simplicity, we identified $E^{11}$ with the a subset of $S^8\times S^7$ by forgetting the first column of $Q$ in \eqref{eq total pullback}. Soon we will see that the quotient of $E^{11}$ by this action is diffeomorphic to the only exotic sphere in dimension 8.

Let us now consider $S^{10}$ as the unit sphere in $\rm{Im}\times \bb H^2$ and the map $f_{10}:S^{10}\to S^7$ defined by
\begin{gather}
f_{10}\begin{pmatrix}p\\w\\x\end{pmatrix}=\begin{pmatrix}\mathbf{b}(p,w)\\x\end{pmatrix}
\end{gather}
where $\mathbf{b}:D^7\to D^4$ is the radial extension of the Blakers-Massey element $b:S^6\to S^3$ firstly defined in \cite{DMR} by
\begin{gather*}
b(p,w) = 
\begin{cases}  
 \frac {w}{|w|}e^{\pi p} \frac{\bar w}{|w|},&  w \neq 0\\
  -1 & w = 0 \, .
\end{cases}
\end{gather*} 
We endow $S^{10}$ with the $S^3$ action
\begin{gather}
q\cdot\begin{pmatrix}p\\w\\x\end{pmatrix}=\begin{pmatrix}p\\qw\\qx\bar q\end{pmatrix}
\end{gather}
and observe that $f_{10}$ is equivariant with respect to this action and \eqref{eq GM action}. By replacing $\mathbf{b}$ by a equivariant smoothing of the mentioned radial extension we get a new bundle $\pi^{13}:f_{10}^*Sp(2)=E^{13}\to S^{10}$ which admits a $\star$-action analogous to the one in $E^{11}$. This action will be also free and we will prove here that

\begin{theorem}\label{thm 8 10}
The quotients of $E^{11}$ and $E^{13}$ by their respective $\star$-actions are diffeomorphic to the only exotic sphere of dimension 8 and a generator of the index 2 subgroup of 10-dimensional homotopy spheres. 
\end{theorem}

By using the notion a construction related to these pullbacks, we also get explicit representatives of their clutching diffeomorphisms. With this in hand we also present an explicit isotopy from the square of the diffeomorphism representing the 8-dimensional sphere to the identity.

\section{Recognizing some Exotic Spheres}


An advantage of describing a manifold as a quotient is a relative freedom about charts. In fact, we get different  descriptions of a quotient by choosing different transversal sections on the total space. Since the total spaces we deal with are, in some sense, pullbacks of other well known spaces, we can try to pullback some good descriptions of one space to the other.

In fact, that is how we can relate the following constructions of our spheres which appears in literature: 
\begin{itemize}
\item the gluing of two discs by an `exotic' diffeomorphism;
\item its plumbing description;
\item the quotient of a relatively reasonable space (e.g. a sphere bundle over a sphere) by an action.
\end{itemize} 

In this section, we use the idea of pullback to induce a plumbing description on $\Sigma^8$ from a plumbing description on $\Sigma^7$. In the case of the Gromoll-Meyer sphere, the passage from quotient to plumbing description is essentially Theorem 1 of \cite{gm}.



We start by choosing suitable subsets on $S^7$ which describe, in a equivariant way, the bundle $\pi:Sp(2)\to S^7$ respecting action $\star$. Here we observe that we can give a full description of this bundle and, therefore, its pullbacks as well. After that, we see how this provides a plumbing description for the quotients.

We can produce a 7-sphere by identifying the common boundary of $D^4\times S^3$ with $S^3\times D^4$ and this is done in the following manner: consider the following subsets of $S^7$:
\begin{align}
\nu (S^3)_+&=\{(x,y)\in S^7~|~x\neq 0\}\\
\nu (S^3)_-&=\{(x,y)\in S^7~|~y\neq 0\}
\end{align}

This set intersect on a tubular neighborhood of the $S^3\times S^3$ torus inside $S^7$ given by $|x|^2=|y|^2=1/2$. So $S^7$ is identified with the gluing of the subsets above along this torus. Topologically this can be also done discarding part of the intersection of each subset $\nu(S^3)_{\pm}$ with the tubular neighborhood of the torus, in a way that the only points which appear on both parts of the manifold is the points of the torus. To say, we can consider, instead of $\nu(S^3)_\pm$, the subsets
\begin{align}
\bar \nu (S^3)_+&=\{(x,y)\in S^7~|~|x|^2\geq 1/2\}\\
\bar \nu (S^3)_-&=\{(x,y)\in S^7~|~|y|^2\geq 1/2\}
\end{align}
Since these subsets are diffeomorphic to $S^3\times D^4$ and $D^4\times S^3$ (being the last the product of a sphere with a disc with radius 1), we can just consider $S^7$ as the gluing of $D^4\times S^3$ and  $S^3\times D^4$ by there common boundary, simultaneously identified with the torus $S^3\times S^3$ thus making $S^7$ the `trivial' plumbing.

A strong point here that makes our calculations considerably cleaner is that we can really make the last identification in the differentiable category. In fact, there is only one way to make this gluing a differentiable manifold preserving the already fixed structures on the products of spheres and discs and the isotopic class of the gluing diffeomorphism of $S^3\times S^3$ (which is, in this case, the identity). Equivalently, there is only one way to \emph{``straighten the angles''} in the resulting manifold.

Now, trivializing $Sp(2)$ over the subsets $\nu(S^3)_\pm$ is simple and we can do it in a way that the $\star$-action has also a nice form. Spelling out this fact we get

\begin{theorem}\label{thm sp2}
$Sp(2)$ with the $\bullet$ and $\star$ actions is equivariantly diffeomorphic to
\begin{align}
D^4\times S^3\times S^3\cup_{f_\theta}S^3\times D^4\times S^3
\end{align}
where $f_\theta(x,y,g)=(x,y,gx\bar y)$ and, in each part the $\bullet$ and $\star$ actions are realized as
\begin{align}
q\bullet(x,y,g)&=(x,y,qg);\\
q\star(x,y,g)&=(qx\bar q,qy\bar q,g\bar q).
\end{align}
\end{theorem}

Here we can think that we realized $Sp(2)$ as the bundle over $D^4\times S^3\cup S^3\times D^4$ with transition function $\theta:S^3\times S^3\to S^3$ given by $\theta(x,y)=x\bar y$. To get an analogous description of a pullback of this bundle we can observe that if we take the first pair of subspaces $\nu(S^3)_\pm$ and consider its pre-images under $f_8$, we get a pair of trivializing subsets on $f_8^*Sp(2)\to S^8$:
\begin{align}
f_8^{-1}(\nu(S^3)_+)&=\{(x,y)\in S^7~|~\lambda^2+|x|^2\neq 0\}\\
f_8^{-1}(\nu(S^3)_-)&=\{(\lambda,x,y)\in S^8~|~w\neq 0\}
\end{align}
Following the same lines of $S^7$, with these subsets we arrive to description of $S^8$ as the glueing of $S^4\times D^4$ and $D^5\times S^3$. Furthermore, we can make its intersection be the exact pre-image of the torus we fixed on $S^7$. With this procedure we get $E^{11}$ with transition map $\theta f_8:S^4\times S^3\to S^3$ given by
\[\theta f_8(\lambda,x,w)=\eta(\lambda, x)w^{-1}\]
where $\eta:S^4\to S^3$ is the suspension of Hopf given by
\[\eta(\lambda,x)=\frac{1}{\sqrt{\lambda^2+|x|^4}}(\lambda+ xi\bar x)\]
The easiest way to get to this maps is through an equivariant homotopy of $f_8$. Indeed, it is a straightforward exercise to observe that equivariant homotopy does not change the diffeomorphism type of the quotient.

The analogous procedure can be done to $f_{10}$. In this case $\theta f_{10}(p,w,y)=b(p,w)\bar y$. In any case, the transition function behaves in the following manner: let $X$ be a $G$-manifold and let $t:X\to G$ satisfy
\begin{gather}\label{eq eqv}
t(gx)=gt(x)g^{-1}
\end{gather}
This property can be seem as a sole responsible for the existence and form of the $\star$-action. Furthermore, this property allows us to produce an equivariant diffeomorphism out of a \emph{`reentrance procedure'} (\cite{DR}): given $t$ as above, we define 
\begin{align*}
\hat t:X&\to X\\x&\mapsto t(x)x 
\end{align*}

Equivariant homotopies of $t$ produces isotopies of $\hat t$; we will use this fact in theorem \ref{thm 13}.

It happens that the quotient by the $\star$-action is completely described by this diffeomorphism and this is a quite general fact (which does not depend on $S^7$ or $Sp(2)$). In what follows, we consider a more general context which reduces to the case of pullbacks of $Sp(2)$ when $G=S^3$.

Let $G$ acts on the unitary standard spheres $S^k$ and $S^l$ by linear actions and consider $S^n=D^{k+1}\times S^l\cup S^k\times D^{l+1}$ as the $G$-manifold obtained by extending linearly the product action on $S^k\times S^l$. Also let $a:S^k\to G$ and $b:S^l\to G$ be maps satisfying \eqref{eq eqv}. Then, we define $r:S^k\times S^l\to G$ as $r(x,y)=a(x)b(y)^{-1}$. We observe that $r$ is equivariant in the sense of \eqref{eq eqv} with respect to the product action on $S^k\times S^l$ so the map $\hat r$ still makes sense.

We now consider an equivariant bundle over $S^n$ using $r$ as clutching map as follows. Let $P_r$ be defined as:
\begin{gather}
P_r=D^{k+1}\times S^l\times G\cup_{f_r}S^k\times D^{l+1}\times G
\end{gather}
where $f_r(x,y,g)=(x,y,gr(x,y))$. We see that the projection $(x,y,g)\mapsto (x,y)$ clearly defines a principal fibration $P_r\to S^{n}$ and define a new action on $P_r$ via the local expression
\begin{equation}\label{eq action star}
q\star(x,y,g)=(qx,qy,gq^{-1})
\end{equation}
This defines a global action since $f_r(q\star(x,y,g))=q\star f_r(x,y,g)$. This action is free, since it is free in the $G$ coordinate and an interesting point of this construction is that we can, somehow, describe the quotient:

\begin{theorem}\label{thm star quotient}
The quotient of $P_r$ by \eqref{eq action star} is diffeomorphic to 
\[\Sigma_r=D^{k+1}\times S^l\cup_{\hat r} S^k\times D^{l+1}\]
where $\hat r:S^k\times S^l\to S^k\times S^l$ is defined by $\hat r(x,y)=(r(x,y)\cdot x,r(x,y)\cdot y)$.\end{theorem}
\begin{proof}
See \cite{DRS}.
\end{proof}

Grouping action \eqref{eq action star} and the principal fibration $P_r\to S^n$, we end with a diagram analogous to \eqref{diagram DCD}.

\begin{equation}\label{diagram rCD}
\begin{xy}\xymatrix{& G\ar@{.}[d]^{\bullet} &\\ G\ar@{..}[r]^{\star} &P_r\ar[d]^\pi\ar[r]^{\pi'} &\Sigma_r\\ &S^{n}&}\end{xy}
\end{equation}

The special way on how the Gromoll-Meyer action descends to $S^7$ allow us, when $G=S^3$, to realize these principal bundles as its pullbacks.

In fact, let $a$ and $b$ be as above with $G=S^3$  and define its \emph{join product} $a*b:S^{l+k+1}\to S^7$  in join coordinates by
\begin{align}
a*b[x,y,t]=[a(x),b(y),t]
\end{align}
We observe that this map can be easily smoothed using mollifiers and the suitable identification of $S^n=S^{l+k+1}$ we do below. 


We recall that the join product of two spheres can be defined by the quotient map $q:S^l\times S^k\times[0,\pi/2]\to S^n$:
\begin{gather}\label{eq join coor}
q(x,y,t)=(\cos t x,\sin ty)
\end{gather}

So, if $f=a*b$,   $f^{-1}(D^{4}\times S^3)=D^{l+1}\times S^k$ and $f^{-1}(S^{3}\times D^4)=S^{l}\times D^{k+1}$. In particular, if $G=S^3$, $f^*Sp(2)\to S^n$ is the bundle defined by
\[P=D^{l+1}\times S^k\times S^3\cup S^l\times D^{k+1}\times S^3\]
with transition map $\theta\circ f(x,y)=a(x)b(y)^{-1}=r(x,y)$, thus smoothly identifying $P_r$ with $f^*Sp(2)$. From our choice of trivializing subsets in $P$ as pre-images of trivializing subsets in $Sp(2)$, there is a tautological identification of $P$ with $f^*Sp(2)$ via the embeddings
\begin{align*}
j_1:f^{-1}(D^4\times S^3)\times S^3&\to M\times Sp(2)\\
(x,g)&\mapsto (x,g\bullet i_1(x))\\
j_2:f^{-1}(S^3\times D^4)\times S^3&\to M\times Sp(2)\\
(x,g)&\mapsto (x,g\bullet i_2(x))
\end{align*}
where $i_1$ and $i_2$ are the embeddings defined by the trivializations $D^4\times S^3\times S^3\to Sp(2)$ and $S^3\times D^4\times S^3\to Sp(2)$. Realizing our identification as an equivariant isomorphism of bundles, therefore, completely identifying the two bundles.

It is also not difficult to see that $f_8$ and $f_{10}$ are equivariantly homotopic to $\rm{id}*\eta$ and $b*\rm{id}$. We set as our next task, the identification of the differentiable class of these manifolds. For this, in what follows, we work directly with a general $r$.

\begin{prop}\label{lem plumbing}
If $r:S^k\times S^l\to G$ is as in \ref{thm star quotient}, then
\begin{equation}\nonumber\hat r=(\hat a\times\rm{id})\circ g_b^{-1}f_a\circ (\rm{id}\times \hat b^{-1})\end{equation}
where $g_b$ and $f_a$ are the diffeomorphisms of $S^k\times S^l$ defined by \begin{align}
g_b(x,y)&=(b(y)x,y)\\f_a(x,y)&=(x,a(x)y)
\end{align} 
\end{prop}

Observing that $\hat a\times\rm{id}$ and $\rm{id}\times \hat b^{-1}$ extends to $S^k\times D^{l+1}$ and $D^{k+1}\times S^l$, respectively, we conclude that $\Sigma_r$ is diffeomorphic to the manifold obtained by gluing $D^{k+1}\times S^l$ with $S^k\times D^{l+1}$ with diffeomorphism $g_b^{-1}f_a$. We also observe that, since the $G$ actions on $S^k$ and $S^l$ are linear, they define group homomorphisms $\Delta_1:G\to SO(k+1)$ and $\Delta_2:G\to SO(l+1)$. By recalling the definition of plumbing, we observe directly that the diffeomorphism $g_b^{-1}f_a$ define the plumbing manifold formed by the homotopy classes of $\Delta_1 b:S^l\to SO(k+1)$ and $\Delta_2 S^k\to SO(l+1)$. We state it as a corollary.
\begin{cor}
The manifold produced by $\hat r$ is diffeomorphic to the plumbing defined by the homotopy classes of $\Delta_2 a$ and $\Delta_1 b$.
\end{cor}

\begin{proof}[Proof of Porposition \ref{lem plumbing}:]
We first observe that, for $f_a(x,y)=(x,a(x)y)$, $f_a(\hat a\times\rm{id})=(\hat a\times\rm{id})f_a=\hat A$, where $A(x,y)=a(x)$, and the analogous for $b$ and $B(x,y)=b(y)$. In particular,
\begin{align*} \hat r&=\hat B^{-1}\hat A\\&= g_b^{-1}(\rm{id}\times \hat b^{-1})f_ a(\hat a\times\rm{id})\\&= g_b^{-1}(\hat a\times \hat b^{-1})f_a\end{align*}
However, $g_b^{-1}(\hat a\times\rm{id})=(\hat a\times\rm{id})g_b^{-1}$ and $(\rm{id}\times \hat b^{-1})f_a=f_a(\rm{id}\times \hat b^{-1})$ since, for example in the first term,
\begin{align}\label{proof tau}f_a^{-1}(\rm{id}\times \hat b^{-1})f_a(x,y)=(x,a(x)^{-1}\hat b^{-1}(a(x)y))=(x,\hat b^{-1}(y)),\nonumber \end{align}
from the equivariance of $\hat b$. So
\[g_b^{-1}(\hat a\times \hat b^{-1})f_a=g_b^{-1}(\hat a\times \rm{id})(\rm{id}\times \hat b^{-1})f_a=(\hat a\times \rm{id})\circ g_b^{-1}f_a\circ (\rm{id}\times \hat b^{-1})\qedhere\]

\end{proof}

In particular, in the case of $f_8$, $\Sigma_r$ is diffeomorphic to the plumbings of the suspension of the generator of $\pi_4SO(3)$ and a generator of $\pi_3SO(5)$. For $f_{10}$ we have the plumbing of the suspension of a generator of $\pi_6SO(3)$ and a generator of $\pi_3SO(7)$. According to \cite{sch-circle}, these are the exotic 8 and 10 spheres, respectively.

\section{A Gromoll-Meyer construction for Kervaire Spheres}\label{sec Kervaire}

There is another family of spheres that also admit a  description as quotients of sphere bundles over spheres. This is the well known family of Kervaire spheres which we describe here  as plumbed manifolds following \cite{brebook}. Let $\tau:S^{n-1}\to O(n)$ be the characteristic class of the tangent bundle of the sphere (or, equivalently, its frame bundle). This class has an special representative defined as
\begin{gather}
\tau(x)v=2\lr{x,v}x-v
\end{gather}
We readily see that $\tau$ is equivariant by $O(n)$ in the following sense:
\[\tau(gx)=g\tau(x)g^{-1}\]
From \cite{brebook}, we know that the Kervaire sphere (or, more generaly, the Kervaire manifold) can be expressed as
\begin{gather}
\Sigma^{2n-1}=D^{n}\times S^{n-1}\cup_{g_\tau^{-1}f_\tau}S^{n-1}\times D^{n} 
\end{gather}
 where $f_\tau$ and $g_\tau$ are analogous to the ones defined on proposition \ref{lem plumbing}. We follow a similar approach to realize these manifolds as possibly interesting quotients. In fact, from theorem \ref{thm star quotient} we see that $\Sigma^{2n-1}$ is diffeomorphic to a quotient of the bundle $P_r$ where $r(x,y)=\tau(x)\tau(y)$ and in what follows, we identify $P_r$ as the pullback of the homogeneous bundle $O(n+1)\to S^{n}$.
 
 We recall that the projection $pr_n:O(n+1)\to S^n$ onto the first column induces a $O(n)$-principal bundle structue. The principal action is realized as right multiplication by block diagonal ortogonal matrices whose first block is the $1\times 1$  matrix  1 and the other is any $O(n)$ matrix. 
 
 It also admit an equivariant action induced by left multiplication by the same type of block diagonal matix. We have, for any $O(n)$ matrices $g,h$ realized as such block diagonal:
 \begin{gather}\label{eq action K}pr_n(gQh)=gpr_n(Q)\end{gather}

This induces the standarrd linear $O(n)$ action on $S^{n}$ by fixing a vector $e_0$.
 
Now, to get the right map we should have in mind that the pre-image  of the equator on $S^n$ must be the torus $S^{n-1}\times S^{n-1}$ on $S^{2n-1}$. This is accomplished, for example, with representatives on the image of the $J$-homomorphisms (think for example on constructing such representatives via Thom-Pontrjagyn). Explicitly, we define the continuous map $f_{2n-1}:S^{2n-1}\to S^n$ as
 \begin{gather}
 f_{2n-1}(x,y)=\exp_{e_0}\pi\tau\Big(\frac{y}{|y|}\Big)x
 \end{gather}

This map is equivariant by the standard 2-axial action of $O(n)$ on $S^{2n-1}$ and \eqref{eq action K}. Indeed
\[f_{2n-1}(gx,gy)=\exp_{e_0}\left(\pi\tau\Big(\frac{gy}{|gy|}\Big)gx\right)=\exp_{e_0}\left(\pi g\tau\Big(\frac{y}{|y|}\Big)x\right)=g\exp_{e_0}\pi \tau\Big(\frac{y}{|y|}\Big)x
\]
since $g$ fixes $e_0.$

Spelling out $\exp_{e_0}\pi x=\cos(\pi |x|)e_0+\sin (\pi|x|)x/|x|$, we observe that the preimage of the equator with respect to the $e_0$ coordinate is formed by the vectors $(x,y)\in S^{2n-1}$ with $|x|=1/2$. This is exactly a torus $S^{n-1}\times S^{n-1}$ and induce a decomposition 
\[S^{2n-1}=D^{n}\times S^{n-1}\cup S^{n-1}\times D^{n} \]
Accordingly to all our proceeding discussion up to now, the clutching function of the pullback bundle will be given by $\theta_{2n-1}(x,y)=\tau f_{2n-1}(x,y)$ and the last is, in fact
\begin{align*}
\tau f_{2n-1}(x,y)&=\tau\left( \exp_{e_0}\pi\tau\Big(\frac{y}{|y|}\Big)x \right)=\tau\left(\tau\Big(\frac{y}{|y|}\Big)x \right)\\
&=\tau\Big(\frac{y}{|y|}\Big)\tau(x) \tau\Big(\frac{y}{|y|}\Big)
\end{align*}
which is actually equal to $\tau(y)\tau(x)\tau(y)$ after identification of the torus of unitary spheres with the embedded one (note that $\tau(x)^{-1}=\tau(x)$). The equivariant bundle formed by such map is indeed isomorphic to $P_r$, however, even easier to see is that the induced diffeomorphism on the quotient is 
 \[(\hat \tau\times\rm{id})\circ g_\tau^{-1}f_\tau\circ (\rm{id}\times \hat \tau)\circ (\rm{id}\times \hat \tau)g_\tau\]
 and, by the same arguments following proposition \ref{lem plumbing}, the manifold glued by this diffeomorphism is diffeomorphic to the one glued by $g_\tau f_\tau$. According to \cite{brebook}, section 1.7, this defines the Kervaire manifold of dimension $2n-1$ which is homeomorphic (and usually not diffeomorphic) to the sphere if $n$ is odd.

We observe here that the biaxial action of $O(n)$ on $S^{2n-1}$ has no fixed points. This is an interesting property if we have the application of the next section in mind. However, when $n$ is odd, we can obtain $S^{n}$ as the quotient of $U(k+1)$ by $U(k)$, where $n=2k+1$. If we carry on with an analogous construction, we end with a biaxial action of $U(k)$ on $S^{4k+1}$ which has fixed points thus putting all Kervaire sphere in the settings of the next section.



\section{Wiedersehen Metrics and Exotic Diffeomorphisms}

It is known (\cite{besse2}) that the behavior of the geodesics in the round sphere is of a very special kind. This behavior can be seem in many forms and one is the motivation of the following definition: a point $x_0\in M$ in a Riemannian manifold is called Blaschke if its cut locus is at constant distance from $x_0$ and wiedersehen if, in addition, its cut locus is a point. We call a manifold which admits a Blaschke (wiedersehen) point as a pointed Blaschke (wiedersehen) manifold. The topological structure of such manifold is given by

\begin{theorem}[Weinstein, Allamigeon-Warner]
A manifold $M^n$ is a pointed Blaschke manifold, if and only if $M$ is diffeomorphic to $D\cup_{\sigma} E$ where $D$ is a disc and $E$ is a disc bundle whose boundary is diffeomorphic to $S^{n-1}$ being $\sigma:\partial D^n\to\partial D(\xi)$ a clutching diffeomorphism.
\end{theorem}

In particular, if $M$ is pointed wiedersehen, the disc bundle above is the trivial disc bundle over a point and $\sigma:S^{n-1}\to S^{n-1}$ is a clutching diffeomorphism that defines $M$ as an homotopy sphere.


The construction of a pointed Blaschke metric on such a manifold was carried out by Weinstein (\cite{besse2}) using the existence of the clutching diffeomorphism. In \cite{dur}, the process was reverted. Starting from an explicit Blaschke metric on the Gromoll-Meyer exotic sphere, a formula for its clutching diffeomorphism is written, inheriting symmetries of the Gromoll-Meyer sphere.

In this section, we observe that the method used on \cite{dur} fits in a more general context.  diagram \ref{diagram DCD} works as a translator of `horizontal' geodesics between $M$ and $M'$. In fact, if we provide $P$ with a metric invariant by both $\star$ and $\bullet$ actions, \emph{bihorizontal} geodesics, i.e., geodesics horizontal with respect to both actions, descend to both $M$ and $M'$ as geodesics which are orthogonal to orbits. We prove

\begin{theorem}\label{thm wiedersehen}
Let $M$ and $M'$ be manifolds fitting in the cross-diagram \eqref{diagram DCD}. Then, $M$ has a fixed point $x_0$  which is  Blaschke (wiedersehen) for a $G$-invriant metric if and only if $M'$ does. 
\end{theorem}

This provides a big family of explicit pointed wiedersehen metrics on ours exotic 8, 10 and Kervaire spheres and also all exotic 7-spheres in \cite{DPR}. We follow some of the ideas in \cite{dur,ADPR} and present some formulas for their attaching diffeomorphisms.

An important point in the proof of theorem \ref{thm wiedersehen} is that we can find a special trivialization around $x_0$ from where we readily get some relevant properties.

Let $\pi:P\to M$ be a $G$-equivariant $G$-bundle which fits in a cross-diagram like \ref{diagram DCD}. We first observe that since both actions on $P$ commute, $P$ can be regarded as a $G\times G$ manifold and we can find a connection 1-form $\omega$ on $\pi$ such that, if $m$ is a $G$-invariant metric on $M$, the connection metric induced by $m$ and $\omega$ is invariant by the $G\times G$-action. In fact, we can average any connection 1-form by the $\star$ action and get the desired result. Let us state it as a lemma.
\begin{lem}\label{lem connection}
If $M$ is a $G$-manifold and $\pi:P\to M$ is a $G$-equivariant $G$-principal bundle then there exists a connection 1-form on $\pi$ which is invariant by the non-principal $G$-action. In particular, for any $G$-invariant metric $m$ on $M$, there exists a principal $G\times G$-invariant metric on $P$ such that the submersion metric induced by $\pi$ is identical to $m$. 
\end{lem}

We also observe that if $x_0$ is a fixed point then, for every  $p\in\pi^{-1}(x)$, there exists an automorphism $\phi:G\to G$ such that the isotropy group at $p$ is the graph of $\phi$, i.e., 
\[(G\times G)=\{(r,\phi(r))\in G\times G~|~r\in G\}\]
A crucial observation here is the following.
\begin{lem}\label{lem iso}
If there is an automorphism $\phi:G\to G$ such that $(G\times G)_p$ is the graph of $\phi$, then $(G\times\{\rm{id}\})p=(\{\rm{id}\}\times G)p=(G\times G)p$. I.e., both actions have the same orbit at $p$.
\end{lem}

It is also not difficult to see that $x_0'=\pi'(p)\in M'$ is also a fixed point on $M'$.
So,  for simplicity, we change the $\bullet$-action so that $(G\times G)_p$ is the diagonal of $G$ on $G\times G$, i.e., is the graph of the identity map. Note that any metric invariant by the old actions is still invariant by the new actions. 

Now, fixing a pointed Blaschke (wiedersehen) metric $m$ of $M$ at $x_0$, we consider a $G\times G$-invariant metric on $P$ as in lemma \ref{lem connection} and its horizontal lift $\cal L_p:T_{x_0}M\to T_pP$. Since the orbits of both $\star$ and $\bullet$ actions coincide at $x_0$, the exponential map  $\exp^P_p:T_pP\to P$, restricted to $\cal L_p(T_{x_0}M)$ provides geodesics which are all the time orthogonal to the orbits of both actions, i.e., it provides geodesics horizontal to both fibrations.  In particular 
\[\pi\exp^P_p(v)=\exp^{M}_{x_0}(d\pi v)\qquad \pi'\exp^P_p(v)=\exp^{M'}_{x_0'}(d\pi' v)\]
The most surprising fact is the following.
\begin{theorem}\label{thm embedding}
If $D\subset T_{x_0}M$ is a disc of radius $\rho$ which is embedded on $M$ by $\exp^M_{x_0}$ then $d\pi'\cal L_p(D)\subset T_{x_0'}M'$ is a disc of radius $\rho$ embedded on $M'$ by $\exp^{M'}_{x_0'}$. 
\end{theorem}
In particular, by assuming that $M$ is Blaschke (\emph{wiedersehen}) at $x_0$, $M'$ will be so at $x_0'$. 

To avoid a cumbersome notation we use the following in the proofs
\[r\star p=rp,\qquad g\bullet p=pg^{-1}\]
With this convention, for $p$ with the chosen isotropy group, $gpg^{-1}=p$, for every $g\in G$.
\begin{proof}[Proof of Theorem \ref{thm embedding}:] Observe that $\exp_p^P:T_pP\to P$ is a $G\times G$-equivariant map since we chose an invariant metric on $P$. The composition of this maps with $\cal L_p$ together with the principal action defines a map $\Psi:T_{x_0}M\times G\to P$ by
\[\Psi(x,g)=\exp_p^P(\cal L_px)g^{-1}=\exp_{pg^{-1}}^P(\cal L_{pg^{-1}}x)\]
Since $\pi(\Psi(x,g))=x$, if we restrict it to a subset $D\subset T_{x_0}M$, $\Psi$ will be a surjective diffeomorphism onto its image as far as the restriction $\exp^M_{x_0}|_D$ is. In particular, if $D$ is a disc of radius $\rho$ which is embedded through $\exp_{x_0}^M$ then $\Psi|_{D\times G}$ is an embedding. To see that  $\pi'(\Psi(D\times G))$ is an embedded disc as well, we actually identify the quotient of $\Psi(D\times G)$ by the $\star$-action with $D$. Once this is done, it is an easy exercise to check that this new disc embedded on $M'$ has radius $\rho$ and its complement is exactly the cut-locus of $x_0'$. We leave this part to the reader. (note that $\overline{\exp_{x_0}(D)}=M'$).
We claim that.
\begin{lem}
$\Psi(rx,gr^{-1})=r\star \Psi(x,g)$
\end{lem}
\begin{proof}
It is straightforward:
\begin{align*}
\Psi(rx,gr^{-1})&=\exp^P_p(\cal L_prx)rg^{-1}
=r\exp^P_{r^{-1}p}(\cal L_{r^{-1}p}(x))rg^{-1}\\
&=r\exp^P_{pr^{-1}}(\cal L_{pr^{-1}}(x))rg^{-1}
=r\exp^P_{p}(\cal L_{p}(x))r^{-1}rg^{-1}\\
&=r\exp^P_{p}(\cal L_{p}(x))g^{-1}=r\Psi(x,g)
\end{align*}
where the second equality holds from the equivariance of $\exp^P_p$ with relation to the $\star$ action.
\end{proof}
In particular we get that the quotient of $\Psi(D\times G)$ by action $\star$ is diffeomorphic to the quotient of $D\times G$ by $r\star(x,g)=(rx,gr^{-1})$. This last quotient is the known usual action used to construct an associated bundle and has the same quotient as the original one, i.e., in our case, it is diffeomorphic to $D$.

\end{proof}

The main interest of theorem \ref{thm embedding} is to provide a topological construction of pointed Blaschke (wiedeserhen) metrics on the $M'$ manifolds. This is a sort of \emph{exoctification} of a original, well-understood, manifold.

We now shift our attention to describe exotic diffeomorphisms of the spheres produced by pullbacks. This is actually a simpler task and follows directly from

\begin{theorem}[\cite{dur}]
The bundle $pr:Sp(2)\to S^7$ is isomorphic to 
\[D^7\times S^3\cup_{f_b} D^7\times S^3\to S^7\]
where $f_b(x,g)=(x,gb(x))$ and $b:S^6\to S^3$ is the Blackers-Massey element defined by
\begin{gather*}
b(p,w) = 
\begin{cases}  
 \frac {w}{|w|}e^{\pi p} \frac{\bar w}{|w|},&  w \neq 0\\
  -1 & w = 0 \, .
\end{cases}
\end{gather*} 
where $p\in\rm{Im}\bb H$, $w\in\bb H$ and $|p|^2+|w|^2=1$. Furthermore, the $\star$-action is given in each $D^7\times S^3$ as
\[r\star((x,y),g)=((rx\bar r,ry\bar r),gr^{-1})\]
\end{theorem}

Geometrically, we can realize the gluing $D^7\cup D^7=S^7$ by taking the two discs as the upper and lower hemispheres of the unitary sphere $S^7\subset \bb H^2$ with respect to the real part of the first coordinate. So the gluing is done along the subsphere whose first coordinate has vanishing real part, making sense with the definition of $b:S^6\to S^3$.

By doing so, we also realize hemisphere decompositions on both $S^8$ and $S^{10}$:
\begin{align*}
\cal S^7=f_8^{-1}(S^6)&=\{(\lambda,x,w)\in S^8~|~\lambda=0\}\\
\cal S^9=f_{10}^{-1}(S^6)&=\{(x,p,w)\in S^{10}~|~\Re x=0\}
\end{align*}
Noticing that the relatives upper and bottom hemispheres are all invariant submanifolds with respect to the action, we realize the bundles $f_8^*Sp(2)$ and $f_{10}^*Sp(2)$ as the 
\begin{align*}
f_n^*Sp(2)=D^n\times S^3\cup D^n\times S^3
\end{align*}
with transition function $f_{bf_n}(x,g)=(x,gb(f_n(x)))$ for $n=8,10$.

Explicitly we have

\begin{theorem}\label{thm 13}
Let $\theta_8:\cal S^7\to S^3$ and $\theta_{10}:\cal S^{9}\to S^3$ be defined as
\begin{align*}
\theta_8(\lambda,x,w) & = 
\begin{cases}  
 \frac{w}{w}e^{{\pi xi \bar x}} \frac{\bar w}{w},&  w \neq 0\\
  -1 & w = 0 \, .
\end{cases},\\\qquad 
\theta_{10}(x,p,w) & = 
\begin{cases}  
b(p,w)e^{{\pi x}} b(p,w)^{-1},&  w \neq 0\\
  -1 & w = 0 \, .
\end{cases}
\end{align*} 
Then $\hat \theta_8:\cal S^7\to\cal S^7$ and $\hat{\theta}_{10}:\cal S^9\to \cal S^9$ are diffeomorphisms which are not isotopic to the identity and represents the exotic 8-sphere and a generator of the subgroup of 10-spheres that bounds paralellizable manifolds.
\end{theorem}

An application of this theorem is an explicit isotopy from twice of the exotic diffeomorphism of the 7-sphere to the identity. This is produced by the means of an equivariant homotopy from $\theta_8$ to $\theta_8^{-1}$.

In fact, let $\alpha(t)=\cos t i+\sin t j$, where $i$ and $j$ are the usual unitary quaternions. We see that $\alpha(t)\rm{Im}\bb H$ for every $t$ and that $\alpha(\pi)=-\alpha(0)$. The desired isotopy is induced by the extension by $-1$ of
\[\theta_{8,t}(x,w)=
 \frac{w}{w}e^{{\pi x\alpha(t) \bar x}} \frac{\bar w}{w}, 
\]
We have,
\[\hat{\theta}_{8,t}\hat{\theta}_8=\widehat{\theta_8\theta_{8,t}}\]
which is equal to $\hat\theta_8^2$ when $t=0$ and to the identity when $t=\pi$.

\section*{Acknowledgements}
Part of this work in the author's Ph. D. thesis under the supervision of A. Rigas and C. Durán to whom the author would like to thank.

\bibliographystyle{alpha}
\bibliography{bib1}

\begin{thebibliography}{ADPR07}

\bibitem[ADPR07]{ADPR}
U.~Abresh, C.~Durán, T.~Püttmann, and A.~Rigas.
\newblock Wiedersehen metrics and exotic involutions of euclidean spheres.
\newblock {\em Journal f\"ur die Reine und Angewandte Mathematik. Crelle's
  Journal}, 605:1--21, 2007.

\bibitem[Bes78]{besse2}
A.~Besse.
\newblock {\em Manifolds all of whose geodesics are closed}.
\newblock Springer, 1978.

\bibitem[Bre72]{brebook}
G.~Bredon.
\newblock {\em Introduction to Compact Transformation Groups}.
\newblock Academic press, 1972.

\bibitem[DPR10]{DPR}
C.~Durán, T.~Püttmann, and A.~Rigas.
\newblock An infinite family of gromoll-meyer spheres.
\newblock {\em Archiv der Mathematik}, 95:269--282, 2010.

\bibitem[DR09]{DR}
C.~Durán and A.~Rigas.
\newblock Equivariant homotopy and deformations of diffeomorphisms.
\newblock {\em Differential Geometry and Its Applications}, 27:206--211, 2009.

\bibitem[DRM04]{DMR}
C.~Durán, A.~Rigas, and A.~Mendoza.
\newblock Blakers-{M}assey elements and exotic diffeomorphisms of ${S}^6$ and
  ${S}^{14}$ via geodesics.
\newblock {\em Transactions of the American Mathematical Society}, pages
  5025--5043, 2004.

\bibitem[DRS10]{DRS}
C.~Durán, A.~Rigas, and L.~Sperança.
\newblock Bootstrapping {A}d-equivariant maps, diffeomorphisms and involutions.
\newblock {\em Matem\'atica Contempor\^anea}, 35:27--39, 2010.

\bibitem[Dur01]{dur}
C.~Dur\'an.
\newblock Pointed wiedersehen metrics on exotic spheres and diffeomorphisms of
  ${S}^6$.
\newblock {\em Geometriae Dedicata}, pages 199--210, 2001.

\bibitem[GM72]{gm}
D.~Gromoll and W.~Meyer.
\newblock An exotic sphere with nonnegative curvature.
\newblock {\em Annals of Mathematics}, 96:413--443, 1972.

\bibitem[Mil56]{mi}
J.~Milnor.
\newblock On manifolds homeomorphic to the 7-sphere.
\newblock {\em Annals of Mathematics}, 64:399--405, 1956.

\bibitem[Sch72]{sch-circle}
R.~Schultz.
\newblock Circle action on homotopy spheres bounding plumbing manifolds.
\newblock {\em Proc. Amer. Soc.}, 36:297--300, 1972.

\end{thebibliography}
\end{document}